\documentclass{sigcomm-alternate}
\usepackage{amsfonts,latexsym,amsmath,amstext,amssymb,verbatim,epsfig,psfrag}
\usepackage{colordvi,pstcol}
\usepackage{graphicx,psboxit}
\usepackage{psfrag}

\def\Z{\mbox{Z\hspace{-.3em}Z}}

\newtheorem{remark}{Remark}
\newtheorem{theorem}{Theorem}

\begin{document}
\title{Understanding differential equations through diffusion point of view: non-symmetric discrete equations}

\numberofauthors{1}
\author{
   \alignauthor Dohy Hong\\
   \affaddr{Alcatel-Lucent Bell Labs}\\
   \affaddr{Route de Villejust}\\
   \affaddr{91620 Nozay, France}\\
   \email{\normalsize dohy.hong@alcatel-lucent.com}
}

\date{\today}
\maketitle

\begin{abstract}
In this paper, we propose a new adaptation of the D-iteration algorithm to numerically solve the differential equations. This problem can be reinterpreted in 2D or 3D (or higher dimensions) as a limit of a diffusion process where the boundary or initial conditions are replaced by fluid catalysts. It has been shown that pre-computing the diffusion process for an elementary catalyst case as a fundamental block of a class of differential equations, the computation efficiency can be greatly improved. Here, we explain how the diffusion point of view can be applied to decompose the fluid diffusion process per direction and how to handle non-symmetric discrete equations.
The method can be applied on the class of problems that can be addressed by the Gauss-Seidel iteration, based on the linear approximation of the differential equations.
\end{abstract}
\begin{psfrags}
\section{Introduction}\label{sec:intro}
The iterative methods to solve differential equations based on the linear approximation are very well
studied approaches \cite{Johnson_1987}, \cite{Ascher:1998:CMO:551054}, \cite{Podlubny:395913}, \cite{Gear:1971:NIV:540426},
\cite{Smith_1985}, \cite{Golub1996}, \cite{Saad}.
The approach we propose here (D-iteration) is a new approach initially applied to numerically solve the eigenvector
of the PageRank type equation \cite{dohy}, \cite{d-algo}, \cite{dist-test}, \cite{distributed}, \cite{partition}, 
\cite{revisit}.

The D-iteration, as diffusion based iteration, is an iteration method that can be understood as a column-vector
based iteration as opposed to a row-vector based approach. Jacobi and Gauss-Seidel iterations are good examples of
row-vector based iteration schemes. While our approach can be associated to the {\em diffusion} vision, the
existing ones can be associated to the {\em collection} vision.

In this paper, we are interested in the numerical solution for linear equation:
\begin{eqnarray}\label{eq:le}
X &=& P.X + B
\end{eqnarray}
where $P$ and $B$ are the matrix and vector associated to the linear approximation of
differential equations with initial conditions or boundary conditions. 

In \cite{diff}, it has been shown how simple adaptations can make the diffusion
approach an interesting candidate as an alternative iterative scheme to numerically solve
differential equations. 
In \cite{diff-germ}, a new approach based on the pre-computation
of the elementary diffusion limit has been proposed. 
In this paper, we study the case of non-symmetric (from the diffusion point of view)
linear equation and how we can very simply decompose the iteration method per direction
for an improved convergence speed.

\section{General equation in 2D}\label{sec:eq}

We consider the general linear (affine) equation of the form:
\begin{align}\label{eq:DD2D}
&U(n,m) = \alpha(+1,0) U(n-1,m)+\alpha(0,+1) U(n,m-1)\\
&+\alpha(-1,0) U(n+1,m)+ \alpha(0,-1) U(n,m+1) + f(n,m)
\end{align}
for $(n,m)\in\Omega$ and with boundary condition:
$H(n,m) = g(n,m)$ for $(n,m)\in\partial\Omega$. We require that
$\partial\Omega$ includes at least all the boundary positions (frontier)
of $\Omega$.

Let's call this problem DD2D (discrete differential 2D) equation problem.
We define: $|\alpha| = |\alpha(-1,0)|+|\alpha(0,-1)|+|\alpha(+1,0)|+|\alpha(0,+1)|$.

Note that the approach proposed here can be directly extended to higher dimension
and also for neighbour positions in Equation \eqref{eq:DD2D} ($(n-1,m)$ etc)
that are more general (it is just required that they are regular on the domain
where we iterate the equations).

\begin{theorem}
[Stability condition]
If $|\alpha|\le 1$, then DD2D has a unique solution. The solution can be obtained
from the iteration of Equation \eqref{eq:DD2D}.
\end{theorem}
\proof
The proof is straightforward noticing that the matrix associated to DD2D has a
spectral radius strictly less than 1.

Now, we associate to DD2D, the diffusion approach: the associated diffusion approach
consists in (iteratively) applying the elementary diffusion operation on $(n,m)$ by:
\begin{align}\label{eq:DD2D-diff}
H(n,m) &+= F(n,m),\\
F(n+1,m) &+= \alpha(+1,0) F(n,m),\\
F(n,m+1) &+= \alpha(0,+1) F(n,m),\\
F(n-1,m) &+= \alpha(-1,0) F(n,m),\\
F(n,m-1) &+= \alpha(0,-1) F(n,m),\\
F(n,m)   &= 0.
\end{align}

\begin{figure}[htbp]
\centering
\includegraphics[width=5cm]{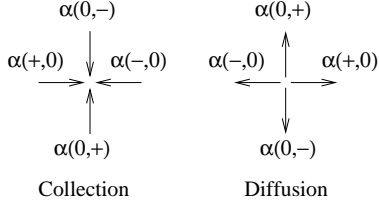}
\caption{Collection/Diffusion approach.}
\label{fig:metal}
\end{figure}

We first build the initial condition for the diffusion process by
applying the elementary diffusion operation to all boundary position in $\partial\Omega$.

The DD2D can be equivalently written under the form:
\begin{eqnarray*}
X &=& P.X + B
\end{eqnarray*}
where the initial fluid $B$ on all positions of $\Omega$ is
defined by the superposition (sum) of $f$ and those coming from the boundary positions.

We recall that the diffusion approach (D-iteration) consists in iteratively
solve the above equation with the fluid vector $F$ and the history vector $H$
(cf. \cite{d-algo}), using \eqref{eq:DD2D-diff}.

\begin{theorem}
Under the stability condition,
the limit of the diffusion scheme defined above converges to the solution of DD2D.
\end{theorem}
\proof
It has been shown in \cite{revisit} that our history vector $H$ corresponds exactly
to $U$ in \eqref{eq:DD2D} when starting from $H=(0,..,0)$ and when applied for the same sequence
of vector entries (than the diffusion process on the fluid vector $F$ \eqref{eq:DD2D-diff},
one vector entry corresponding here to a spatial position $(n,m)$).

Below, we show through simple examples how our approach works very simply and how it can
help us decomposing the diffusion process for a better efficiency.

\subsection{Catalyst limit in 1D}\label{sec:cata1d}
The general case of the linear DD1D operator associated to the diffusion is:
\begin{align}\label{eq:DD1D-diff}
H(n) &+= F(n),\\
F(n+1) &+= \alpha(+) F(n),\\
F(n-1) &+= \alpha(-) F(n),\\
F(n)   &= 0.
\end{align}
Its elementary catalyst limit (cf. \cite{diff-germ}) $\phi$ is associated to
the solution of \eqref{eq:DD1D-diff} with
the initial condition $g(0) = 1$ and with the constraint that the diffusion
is applied once at position $0$ and then the position behaves as a black hole
(diffusion only on $n\neq 0$ and with $0$ at $+\infty$).
The solution is here simple and explicit (for instance for $\alpha(+) > 0$ and
$\alpha(-) > 0$ and $|\alpha|<1$): we can solve the equation
$$
\alpha(-)x^2 - x + \alpha(+) = 0
$$
for the solution on $\Z^+$, which is $\phi_+^{\infty}(n) = r_+^n$ with 
$$r_+ = \frac{1-\sqrt{1-4\alpha(+)\alpha(-)}}{2\alpha(-)}$$
and for $\Z^-$, $\phi_-^{\infty}(-n) = r_-^{n}$ with 
$$r_- = \frac{1-\sqrt{1-4\alpha(+)\alpha(-)}}{2\alpha(+)}.$$

One can easily check that the above solution is the right one for all
cases of $|\alpha|<1$.

If we put a finite boundary condition $g(N)=0$, we find the exact limit
by successive compensations of the surplus or deficit of fluid at the two
boundary positions:
\begin{align*}
\phi_+^N(n) &=  r_+^n - r_+^N r_-^N r_-^{-n} + r_+^N r_-^N r_+^{n}\\
          & - r_+^{2N} r_-^{2N} r_-^{-n} + r_+^{2N} r_-^{2N} r_+^{n} + ...\\
 & = r_+^n\frac{1-(r_+r_-)^{N-n}}{1-(r_+r_-)^{N}},
\end{align*}
and for $\Z^-$,
\begin{align*}
\phi_-^N(-n) &= r_-^n\frac{1-(r_+r_-)^{N-n}}{1-(r_+r_-)^{N}}.
\end{align*}

In particular, if $r=r_+=r_-$, we have on $[-N,N]$:
\begin{align*}
\phi^N(n) &= r^{|n|}\frac{1-r^{2(N-|n|)}}{1-r^{2N}}.
\end{align*}

Note that $\phi^N$ is the elementary function we can use when injecting fluid from the
boundary position. When the value at the origin is not imposed ($0$ is not a black hole), we have
to use the normalized version of $\phi$;
$$
\tilde{\phi}(n) = \frac{\phi(n)}{1-\alpha(-)\phi(1)-\alpha(+)\phi(-1)}
$$ 
on $\Z$ without bound and on $[-N,N]$,
$$
\tilde{\phi^N}(n) = \frac{\phi^N(n)}{1-\alpha(-)\phi^N(1)-\alpha(+)\phi^N(-1)}.
$$ 

\begin{remark}
$\phi(0)=1$ by definition. $\tilde{\phi}(0) \ge 1$ represents the total fluid that comes back to $0$.
$\alpha(-)\phi(1)$ represents the fluid that goes to the black hole coming from $\Z^+$.
$\alpha(+)\phi(-1)$ represents the fluid that goes to the black hole coming from $\Z^-$.
\end{remark}
\subsection{Example of 1D differential equation with time dimension}
Let's consider a concrete case (cf. \cite{webedp}) of heat equation evolution
in time in 1D:
$$
\partial_t U(t,x) = {\partial_x}^2 U(t,x), (t,x)\in [0,T]\times [0,1]
$$
with initial condition $U(O,x) = U_0(x) = sin(\pi x)$ (pre-heated metal stick) and boundary condition
$U(t,0)=U(t,1) = 0$ (imposed temperature at the extreme points of the metal).
Then we can discretize the above equation as (usually called implicit equation):
\begin{align}\label{eq:1DT}
\frac{U(t,n)-U(t-1,n)}{\Delta t} =& \frac{U(t,n+1)+U(t,n-1)-2 U(t,n)}{\Delta x^2}
\end{align}
which can be written as:
\begin{align*}
U(t,n) =& \frac{1}{1+2 k} U(t-1,n) + \frac{k}{1+2k} (U(t,n+1)+U(t,n-1)).
\end{align*}

First it is well known that the above scheme is always stable (cf. \cite{Lax}).
Here, we have $\alpha(-1,0) = 0$, $\alpha(0,+1) = \alpha(0,-1) = \frac{k}{1+2k}$ and
$\alpha(+1,0) = \frac{1}{1+2 k}$.
The initial condition is build by injecting $U_0(x)$ to $F(1,x)$:
\begin{align}\label{eq:pred}
F(1,x) &:= \alpha(+1,0) U_0(x). 
\end{align}

Now, thanks to the freedom of the order in which we apply the diffusion (on the position
choice), we do the following:
\begin{itemize}
\item we first advance on the time axis once using only the diffusion with $\alpha(+1,0)$
  (this is the application of \eqref{eq:pred} for time $t=1$);
\item then we freeze the diffusion on time axis and diffuse only on x-axis;
  since we know the exact limit of the elementary catalyst solution (diffusion
  of 1 from $x=0$ with boundary condition equal to 0 at $N$, which
  is $\tilde{\phi}^N(n)$ (cf. Section \ref{sec:cata1d}) with $r = (1-\sqrt{1-4\alpha^2} )/(2\alpha)$ ($\alpha = \alpha(0,+1)$)
  we can diffuse all fluid for $x=1$ to $x=L_x-1$ using the pre-computed elementary solution, 
  and obtain the exact directional (x-axis) diffusion limit (without the need to compensate the surplus fluid
  at the boundary $x=0$, $x=L_x$); note that we could also use $\tilde{\phi}^{\infty}(n) = r^n/(1-2rk/(1+2k))$, then decide to compensate the
  boundary values using iteratively ${\phi}^{\infty}(n)=r^n$ or using once ${\phi}^{N}(n)$); 
\item then we restart the process for the next time $t+1$ by diffusing $H(t,)$ with $\alpha(+1,0)$: $$F(t+1,x) := \alpha(+1,0) H(t,x).$$ 
\end{itemize}

Note that there is no approximation in the above scheme.

We show the comparison of our approach (DI) to the naive iteration of Equation \eqref{eq:1DT} (GS as
Gauss-Seidel).

\begin{figure}[htbp]
\centering
\includegraphics[angle=-90,width=7cm]{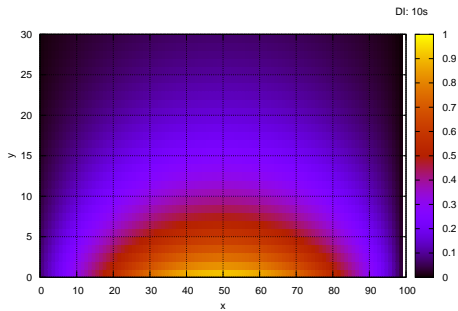}
\caption{D-iteration: 10s (exact limit).}
\label{fig:1DT-DI}
\end{figure}
\begin{figure}[htbp]
\centering
\includegraphics[angle=-90,width=7cm]{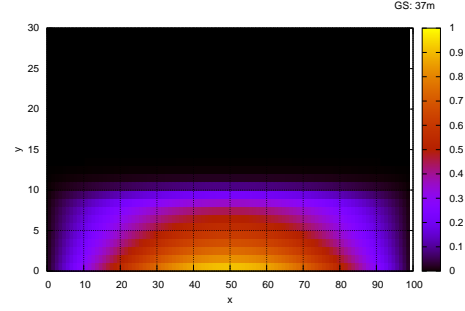}
\caption{Gauss-Seidel: 37m.}
\label{fig:1DT-GS1}
\end{figure}
\begin{figure}[htbp]
\centering
\includegraphics[angle=-90,width=7cm]{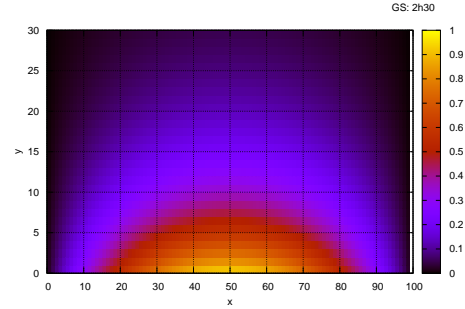}
\caption{Gauss-Seidel: 2h30m. Not yet converged.}
\label{fig:1DT-GS}
\end{figure}

\subsection{Example of 1D differential equation of order 2}
Let's consider the general second order linear differential equation:
$$
y''(x) + \alpha y'(x) + \beta y(x) = f(x).
$$
Its discretized equation is:
\begin{align}\label{eq:1D2}
&\frac{y(n+1)+y(n-1)-2y(n)}{\Delta x^2} + \alpha \frac{y(n)-y(n-1)}{\Delta x}\\
& + \beta y(n) = f(n).
\end{align}
which can be written as:
\begin{align}\label{eq:1D2d}
y(n) =&\alpha(-) y(n+1)+\alpha(+)y(n-1) - \gamma f(n),
\end{align}
with $\alpha(-)=\frac{1}{2+\alpha\Delta x - \beta\Delta x^2}$,
$\alpha(+)=\frac{-\alpha\Delta x}{2+\alpha\Delta x - \beta\Delta x^2}$
and $\gamma = \frac{\Delta x^2}{2+\alpha\Delta x - \beta\Delta x^2}$.

A sufficient stability condition for this equation is: $\Delta \le \min(\frac{2|\alpha|}{|\beta|},1/|\alpha|,1/\sqrt{|\beta|})$.
We can then apply the method as in the previous section using the elementary catalyst
limit associated to $\alpha(+)$ and $\alpha(-)$ (cf. Section \ref{sec:cata1d}).

\section{Conclusion}\label{sec:conclusion}
In this paper we showed that thanks to the diffusion point of view
we can efficiently solve the linear equations associated to
non-symmetric operators and that we could also exploit the idea
of diffusion per direction for a faster computation.
This last idea of diffusion per direction is a promising approach in the
context of linear operators in higher dimension when the naive iteration method
becomes really painful.
Further exploitation of this will be addressed in a future paper.

\end{psfrags}
\bibliographystyle{abbrv}
\bibliography{sigproc}

\end{document}